\documentclass[12pt]{amsart}
\usepackage{amsmath,amssymb,amsfonts}

    \usepackage[english]{babel}
    
    \usepackage[T1]{fontenc}
    \usepackage{mathpazo}

    \usepackage{graphicx}
    \makeatletter
    \def\maxwidth{\ifdim\Gin@nat@width>\linewidth\linewidth
    \else\Gin@nat@width\fi}
    \makeatother
    \let\Oldincludegraphics\includegraphics
    \renewcommand{\includegraphics}[1]{\Oldincludegraphics[width=.8\maxwidth]{#1}}
    \usepackage{caption}
    \DeclareCaptionLabelFormat{nolabel}{}
    \usepackage{adjustbox} % Used to constrain images to a maximum size 
    \usepackage{xcolor} % Allow colors to be defined
    \usepackage{enumerate} % Needed for markdown enumerations to work
    \usepackage{geometry} % Used to adjust the document margins
    \usepackage{amsmath} % Equations
    \usepackage{amssymb} % Equations
    \usepackage{textcomp} % defines textquotesingle
    \AtBeginDocument{%
        \def\PYZsq{\textquotesingle}% Upright quotes in Pygmentized code
    }
    \usepackage{upquote} % Upright quotes for verbatim code
    \usepackage{eurosym} % defines \euro
    \usepackage[mathletters]{ucs} % Extended unicode (utf-8) support
    \usepackage[utf8x]{inputenc} % Allow utf-8 characters in the tex document
    \usepackage{fancyvrb} % verbatim replacement that allows latex
    \usepackage{grffile} % extends the file name processing of package graphics 
    \usepackage{hyperref}
    \usepackage{longtable} % longtable support required by pandoc >1.10
    \usepackage{booktabs}  % table support for pandoc > 1.12.2
    \usepackage[inline]{enumitem} % IRkernel/repr support (it uses the enumerate* environment)
    \usepackage[normalem]{ulem} % ulem is needed to support strikethroughs (\sout)
    \usepackage{mathrsfs}

    \definecolor{urlcolor}{rgb}{0,.145,.698}
    \definecolor{linkcolor}{rgb}{.71,0.21,0.01}
    \definecolor{citecolor}{rgb}{.12,.54,.11}

    \definecolor{ansi-black}{HTML}{3E424D}
    \definecolor{ansi-black-intense}{HTML}{282C36}
    \definecolor{ansi-red}{HTML}{E75C58}
    \definecolor{ansi-red-intense}{HTML}{B22B31}
    \definecolor{ansi-green}{HTML}{00A250}
    \definecolor{ansi-green-intense}{HTML}{007427}
    \definecolor{ansi-yellow}{HTML}{DDB62B}
    \definecolor{ansi-yellow-intense}{HTML}{B27D12}
    \definecolor{ansi-blue}{HTML}{208FFB}
    \definecolor{ansi-blue-intense}{HTML}{0065CA}
    \definecolor{ansi-magenta}{HTML}{D160C4}
    \definecolor{ansi-magenta-intense}{HTML}{A03196}
    \definecolor{ansi-cyan}{HTML}{60C6C8}
    \definecolor{ansi-cyan-intense}{HTML}{258F8F}
    \definecolor{ansi-white}{HTML}{C5C1B4}
    \definecolor{ansi-white-intense}{HTML}{A1A6B2}
    \definecolor{ansi-default-inverse-fg}{HTML}{FFFFFF}
    \definecolor{ansi-default-inverse-bg}{HTML}{000000}

    \providecommand{\tightlist}{%
      \setlength{\itemsep}{0pt}\setlength{\parskip}{0pt}}
    \DefineVerbatimEnvironment{Highlighting}{Verbatim}{commandchars=\\\{\}}
    \let\Oldtex\TeX
    \let\Oldlatex\LaTeX
    \renewcommand{\TeX}{\textrm{\Oldtex}}
    \renewcommand{\LaTeX}{\textrm{\Oldlatex}}
    \title{Experiments with the census}
    \author{Igor Rivin}
    \email{rivin@temple.edu}
    \address{Mathematics Department, Temple University}
    \thanks{This paper was written when the author was invited to present a talk at the CUNY Graduate Center Seminar on Group Theory. The author would like to thank the organizers for their hospitality.}

\makeatletter
\def\PY@reset{\let\PY@it=\relax \let\PY@bf=\relax%
    \let\PY@ul=\relax \let\PY@tc=\relax%
    \let\PY@bc=\relax \let\PY@ff=\relax}
\def\PY@tok#1{\csname PY@tok@#1\endcsname}
\def\PY@toks#1+{\ifx\relax#1\empty\else%
    \PY@tok{#1}\expandafter\PY@toks\fi}
\def\PY@do#1{\PY@bc{\PY@tc{\PY@ul{%
    \PY@it{\PY@bf{\PY@ff{#1}}}}}}}
\def\PY#1#2{\PY@reset\PY@toks#1+\relax+\PY@do{#2}}

\expandafter\def\csname PY@tok@w\endcsname{\def\PY@tc##1{\textcolor[rgb]{0.73,0.73,0.73}{##1}}}
\expandafter\def\csname PY@tok@c\endcsname{\let\PY@it=\textit\def\PY@tc##1{\textcolor[rgb]{0.25,0.50,0.50}{##1}}}
\expandafter\def\csname PY@tok@cp\endcsname{\def\PY@tc##1{\textcolor[rgb]{0.74,0.48,0.00}{##1}}}
\expandafter\def\csname PY@tok@k\endcsname{\let\PY@bf=\textbf\def\PY@tc##1{\textcolor[rgb]{0.00,0.50,0.00}{##1}}}
\expandafter\def\csname PY@tok@kp\endcsname{\def\PY@tc##1{\textcolor[rgb]{0.00,0.50,0.00}{##1}}}
\expandafter\def\csname PY@tok@kt\endcsname{\def\PY@tc##1{\textcolor[rgb]{0.69,0.00,0.25}{##1}}}
\expandafter\def\csname PY@tok@o\endcsname{\def\PY@tc##1{\textcolor[rgb]{0.40,0.40,0.40}{##1}}}
\expandafter\def\csname PY@tok@ow\endcsname{\let\PY@bf=\textbf\def\PY@tc##1{\textcolor[rgb]{0.67,0.13,1.00}{##1}}}
\expandafter\def\csname PY@tok@nb\endcsname{\def\PY@tc##1{\textcolor[rgb]{0.00,0.50,0.00}{##1}}}
\expandafter\def\csname PY@tok@nf\endcsname{\def\PY@tc##1{\textcolor[rgb]{0.00,0.00,1.00}{##1}}}
\expandafter\def\csname PY@tok@nc\endcsname{\let\PY@bf=\textbf\def\PY@tc##1{\textcolor[rgb]{0.00,0.00,1.00}{##1}}}
\expandafter\def\csname PY@tok@nn\endcsname{\let\PY@bf=\textbf\def\PY@tc##1{\textcolor[rgb]{0.00,0.00,1.00}{##1}}}
\expandafter\def\csname PY@tok@ne\endcsname{\let\PY@bf=\textbf\def\PY@tc##1{\textcolor[rgb]{0.82,0.25,0.23}{##1}}}
\expandafter\def\csname PY@tok@nv\endcsname{\def\PY@tc##1{\textcolor[rgb]{0.10,0.09,0.49}{##1}}}
\expandafter\def\csname PY@tok@no\endcsname{\def\PY@tc##1{\textcolor[rgb]{0.53,0.00,0.00}{##1}}}
\expandafter\def\csname PY@tok@nl\endcsname{\def\PY@tc##1{\textcolor[rgb]{0.63,0.63,0.00}{##1}}}
\expandafter\def\csname PY@tok@ni\endcsname{\let\PY@bf=\textbf\def\PY@tc##1{\textcolor[rgb]{0.60,0.60,0.60}{##1}}}
\expandafter\def\csname PY@tok@na\endcsname{\def\PY@tc##1{\textcolor[rgb]{0.49,0.56,0.16}{##1}}}
\expandafter\def\csname PY@tok@nt\endcsname{\let\PY@bf=\textbf\def\PY@tc##1{\textcolor[rgb]{0.00,0.50,0.00}{##1}}}
\expandafter\def\csname PY@tok@nd\endcsname{\def\PY@tc##1{\textcolor[rgb]{0.67,0.13,1.00}{##1}}}
\expandafter\def\csname PY@tok@s\endcsname{\def\PY@tc##1{\textcolor[rgb]{0.73,0.13,0.13}{##1}}}
\expandafter\def\csname PY@tok@sd\endcsname{\let\PY@it=\textit\def\PY@tc##1{\textcolor[rgb]{0.73,0.13,0.13}{##1}}}
\expandafter\def\csname PY@tok@si\endcsname{\let\PY@bf=\textbf\def\PY@tc##1{\textcolor[rgb]{0.73,0.40,0.53}{##1}}}
\expandafter\def\csname PY@tok@se\endcsname{\let\PY@bf=\textbf\def\PY@tc##1{\textcolor[rgb]{0.73,0.40,0.13}{##1}}}
\expandafter\def\csname PY@tok@sr\endcsname{\def\PY@tc##1{\textcolor[rgb]{0.73,0.40,0.53}{##1}}}
\expandafter\def\csname PY@tok@ss\endcsname{\def\PY@tc##1{\textcolor[rgb]{0.10,0.09,0.49}{##1}}}
\expandafter\def\csname PY@tok@sx\endcsname{\def\PY@tc##1{\textcolor[rgb]{0.00,0.50,0.00}{##1}}}
\expandafter\def\csname PY@tok@m\endcsname{\def\PY@tc##1{\textcolor[rgb]{0.40,0.40,0.40}{##1}}}
\expandafter\def\csname PY@tok@gh\endcsname{\let\PY@bf=\textbf\def\PY@tc##1{\textcolor[rgb]{0.00,0.00,0.50}{##1}}}
\expandafter\def\csname PY@tok@gu\endcsname{\let\PY@bf=\textbf\def\PY@tc##1{\textcolor[rgb]{0.50,0.00,0.50}{##1}}}
\expandafter\def\csname PY@tok@gd\endcsname{\def\PY@tc##1{\textcolor[rgb]{0.63,0.00,0.00}{##1}}}
\expandafter\def\csname PY@tok@gi\endcsname{\def\PY@tc##1{\textcolor[rgb]{0.00,0.63,0.00}{##1}}}
\expandafter\def\csname PY@tok@gr\endcsname{\def\PY@tc##1{\textcolor[rgb]{1.00,0.00,0.00}{##1}}}
\expandafter\def\csname PY@tok@ge\endcsname{\let\PY@it=\textit}
\expandafter\def\csname PY@tok@gs\endcsname{\let\PY@bf=\textbf}
\expandafter\def\csname PY@tok@gp\endcsname{\let\PY@bf=\textbf\def\PY@tc##1{\textcolor[rgb]{0.00,0.00,0.50}{##1}}}
\expandafter\def\csname PY@tok@go\endcsname{\def\PY@tc##1{\textcolor[rgb]{0.53,0.53,0.53}{##1}}}
\expandafter\def\csname PY@tok@gt\endcsname{\def\PY@tc##1{\textcolor[rgb]{0.00,0.27,0.87}{##1}}}
\expandafter\def\csname PY@tok@err\endcsname{\def\PY@bc##1{\setlength{\fboxsep}{0pt}\fcolorbox[rgb]{1.00,0.00,0.00}{1,1,1}{\strut ##1}}}
\expandafter\def\csname PY@tok@kc\endcsname{\let\PY@bf=\textbf\def\PY@tc##1{\textcolor[rgb]{0.00,0.50,0.00}{##1}}}
\expandafter\def\csname PY@tok@kd\endcsname{\let\PY@bf=\textbf\def\PY@tc##1{\textcolor[rgb]{0.00,0.50,0.00}{##1}}}
\expandafter\def\csname PY@tok@kn\endcsname{\let\PY@bf=\textbf\def\PY@tc##1{\textcolor[rgb]{0.00,0.50,0.00}{##1}}}
\expandafter\def\csname PY@tok@kr\endcsname{\let\PY@bf=\textbf\def\PY@tc##1{\textcolor[rgb]{0.00,0.50,0.00}{##1}}}
\expandafter\def\csname PY@tok@bp\endcsname{\def\PY@tc##1{\textcolor[rgb]{0.00,0.50,0.00}{##1}}}
\expandafter\def\csname PY@tok@fm\endcsname{\def\PY@tc##1{\textcolor[rgb]{0.00,0.00,1.00}{##1}}}
\expandafter\def\csname PY@tok@vc\endcsname{\def\PY@tc##1{\textcolor[rgb]{0.10,0.09,0.49}{##1}}}
\expandafter\def\csname PY@tok@vg\endcsname{\def\PY@tc##1{\textcolor[rgb]{0.10,0.09,0.49}{##1}}}
\expandafter\def\csname PY@tok@vi\endcsname{\def\PY@tc##1{\textcolor[rgb]{0.10,0.09,0.49}{##1}}}
\expandafter\def\csname PY@tok@vm\endcsname{\def\PY@tc##1{\textcolor[rgb]{0.10,0.09,0.49}{##1}}}
\expandafter\def\csname PY@tok@sa\endcsname{\def\PY@tc##1{\textcolor[rgb]{0.73,0.13,0.13}{##1}}}
\expandafter\def\csname PY@tok@sb\endcsname{\def\PY@tc##1{\textcolor[rgb]{0.73,0.13,0.13}{##1}}}
\expandafter\def\csname PY@tok@sc\endcsname{\def\PY@tc##1{\textcolor[rgb]{0.73,0.13,0.13}{##1}}}
\expandafter\def\csname PY@tok@dl\endcsname{\def\PY@tc##1{\textcolor[rgb]{0.73,0.13,0.13}{##1}}}
\expandafter\def\csname PY@tok@s2\endcsname{\def\PY@tc##1{\textcolor[rgb]{0.73,0.13,0.13}{##1}}}
\expandafter\def\csname PY@tok@sh\endcsname{\def\PY@tc##1{\textcolor[rgb]{0.73,0.13,0.13}{##1}}}
\expandafter\def\csname PY@tok@s1\endcsname{\def\PY@tc##1{\textcolor[rgb]{0.73,0.13,0.13}{##1}}}
\expandafter\def\csname PY@tok@mb\endcsname{\def\PY@tc##1{\textcolor[rgb]{0.40,0.40,0.40}{##1}}}
\expandafter\def\csname PY@tok@mf\endcsname{\def\PY@tc##1{\textcolor[rgb]{0.40,0.40,0.40}{##1}}}
\expandafter\def\csname PY@tok@mh\endcsname{\def\PY@tc##1{\textcolor[rgb]{0.40,0.40,0.40}{##1}}}
\expandafter\def\csname PY@tok@mi\endcsname{\def\PY@tc##1{\textcolor[rgb]{0.40,0.40,0.40}{##1}}}
\expandafter\def\csname PY@tok@il\endcsname{\def\PY@tc##1{\textcolor[rgb]{0.40,0.40,0.40}{##1}}}
\expandafter\def\csname PY@tok@mo\endcsname{\def\PY@tc##1{\textcolor[rgb]{0.40,0.40,0.40}{##1}}}
\expandafter\def\csname PY@tok@ch\endcsname{\let\PY@it=\textit\def\PY@tc##1{\textcolor[rgb]{0.25,0.50,0.50}{##1}}}
\expandafter\def\csname PY@tok@cm\endcsname{\let\PY@it=\textit\def\PY@tc##1{\textcolor[rgb]{0.25,0.50,0.50}{##1}}}
\expandafter\def\csname PY@tok@cpf\endcsname{\let\PY@it=\textit\def\PY@tc##1{\textcolor[rgb]{0.25,0.50,0.50}{##1}}}
\expandafter\def\csname PY@tok@c1\endcsname{\let\PY@it=\textit\def\PY@tc##1{\textcolor[rgb]{0.25,0.50,0.50}{##1}}}
\expandafter\def\csname PY@tok@cs\endcsname{\let\PY@it=\textit\def\PY@tc##1{\textcolor[rgb]{0.25,0.50,0.50}{##1}}}

\def\PYZus{\char`\_}
\def\PYZob{\char`\{}
\def\PYZcb{\char`\}}

\def\PYZgt{\char`\>}

\def\PYZhy{\char`\-}
\def\PYZsq{\char`\'}

\makeatother

    \definecolor{incolor}{rgb}{0.0, 0.0, 0.5}
    \definecolor{outcolor}{rgb}{0.545, 0.0, 0.0}

    \hypersetup{
      breaklinks=true,  % so long urls are correctly broken across lines
      colorlinks=true,
      urlcolor=urlcolor,
      linkcolor=linkcolor,
      citecolor=citecolor,
      }
\begin{document}
    
    \begin{abstract}
        In this paper we study the manifolds in the census of ``small'' 3-manifolds as available in SnapPy. We compare our results with the statistics of random 3-manifolds obtained using the Dunfield Thurston and Rivin models.
    \end{abstract}
    \maketitle

    \hypertarget{Introduction}{%
\section{Introduction}\label{Introduction}}
In past work (\cite{rivin2014statistics}) we have studied the statistics of random manifolds fibering over the circle, using a model (the ``Rivin Model'') similar in spirit to that used by N. Dunfield and W.P.Thurston in \cite{dunfield2006finite}. The distributions of manifolds produced by both models is clearly skewed (the manifolds tend to be "long and skinny".) Many have argued that the "right" model is that of randomly gluing tetrahedra together, then throwing out those gluings that are not manifolds. Unfortunately, this is not at all probabilistically tractable, so we do the next best thing and look at \emph{all} manifolds possessing small triangulations, thanks to the census of such manifolds built into \texttt{SnapPy} - \cite{SnapPy}. This paper started life as a Jupyter notebook.
   \hypertarget{preliminaries}{%
\section{Preliminaries}\label{Preliminaries}}

    First, we import the usual (and some unusual) libraries:

    \begin{Verbatim}[commandchars=\\\{\}]
{\color{incolor}In [{\color{incolor}1}]:} \PY{k+kn}{from} \PY{n+nn}{snappy} \PY{k}{import} \PY{o}{*}
        \PY{k+kn}{from} \PY{n+nn}{multiprocessing} \PY{k}{import} \PY{n}{Pool}
        \PY{k+kn}{import} \PY{n+nn}{pandas} \PY{k}{as} \PY{n+nn}{pd}
        \PY{k+kn}{import} \PY{n+nn}{numpy} \PY{k}{as} \PY{n+nn}{np}
        \PY{k+kn}{import} \PY{n+nn}{functools}
        \PY{k+kn}{from} \PY{n+nn}{operator} \PY{k}{import} \PY{n}{mul}
        \PY{k+kn}{import} \PY{n+nn}{xgboost} \PY{k}{as} \PY{n+nn}{xg}
        \PY{k+kn}{import} \PY{n+nn}{matplotlib}\PY{n+nn}{.}\PY{n+nn}{pyplot} \PY{k}{as} \PY{n+nn}{plt}
        \PY{k+kn}{import} \PY{n+nn}{seaborn} \PY{k}{as} \PY{n+nn}{sns}
        \PY{k+kn}{from} \PY{n+nn}{sklearn}\PY{n+nn}{.}\PY{n+nn}{model\PYZus{}selection} \PY{k}{import} \PY{n}{train\PYZus{}test\PYZus{}split}
        \PY{k+kn}{from} \PY{n+nn}{sklearn}\PY{n+nn}{.}\PY{n+nn}{model\PYZus{}selection} \PY{k}{import} \PY{n}{GridSearchCV}
        \PY{k+kn}{from} \PY{n+nn}{sklearn}\PY{n+nn}{.}\PY{n+nn}{cluster} \PY{k}{import} \PY{n}{KMeans}
        \PY{k+kn}{from} \PY{n+nn}{sklearn}\PY{n+nn}{.}\PY{n+nn}{manifold} \PY{k}{import} \PY{n}{TSNE}
\end{Verbatim}

    Now we define some utility functions to deal with SnapPy's goofy
formats.

    First, the length spectrum

    \begin{Verbatim}[commandchars=\\\{\}]
{\color{incolor}In [{\color{incolor}2}]:} \PY{k}{def} \PY{n+nf}{mung\PYZus{}spec}\PY{p}{(}\PY{n}{specline}\PY{p}{)}\PY{p}{:}
            \PY{n}{mult} \PY{o}{=} \PY{n}{specline}\PY{p}{[}\PY{l+s+s1}{\PYZsq{}}\PY{l+s+s1}{multiplicity}\PY{l+s+s1}{\PYZsq{}}\PY{p}{]}
            \PY{n}{thelen} \PY{o}{=} \PY{n}{specline}\PY{p}{[}\PY{l+s+s1}{\PYZsq{}}\PY{l+s+s1}{length}\PY{l+s+s1}{\PYZsq{}}\PY{p}{]}
            \PY{k}{return} \PY{p}{[}\PY{n}{thelen}\PY{p}{]}\PY{o}{*}\PY{n}{mult}
\end{Verbatim}

    \begin{Verbatim}[commandchars=\\\{\}]
{\color{incolor}In [{\color{incolor}3}]:} \PY{k}{def} \PY{n+nf}{mung\PYZus{}all\PYZus{}spec}\PY{p}{(}\PY{n}{thespec}\PY{p}{)}\PY{p}{:}
            \PY{n}{thelens} \PY{o}{=} \PY{p}{[}\PY{n}{mung\PYZus{}spec}\PY{p}{(}\PY{n}{i}\PY{p}{)}\PY{k}{for} \PY{n}{i} \PY{o+ow}{in} \PY{n}{thespec}\PY{p}{]}
            \PY{k}{return} \PY{n+nb}{sum}\PY{p}{(}\PY{n}{thelens}\PY{p}{,} \PY{p}{[}\PY{p}{]}\PY{p}{)}
\end{Verbatim}

    \begin{Verbatim}[commandchars=\\\{\}]
{\color{incolor}In [{\color{incolor}4}]:} \PY{k}{def} \PY{n+nf}{mung\PYZus{}complexes}\PY{p}{(}\PY{n}{clist}\PY{p}{)}\PY{p}{:}
            \PY{n}{cs} \PY{o}{=} \PY{p}{[}\PY{p}{[}\PY{n+nb}{float}\PY{p}{(}\PY{n}{x}\PY{o}{.}\PY{n}{real}\PY{p}{(}\PY{p}{)}\PY{p}{)}\PY{p}{,} \PY{n+nb}{float}\PY{p}{(}\PY{n}{x}\PY{o}{.}\PY{n}{imag}\PY{p}{(}\PY{p}{)}\PY{p}{)}\PY{p}{]} \PY{k}{for} \PY{n}{x} \PY{o+ow}{in} \PY{n}{clist}\PY{p}{]}
            \PY{k}{return} \PY{n+nb}{sum}\PY{p}{(}\PY{n}{cs}\PY{p}{,} \PY{p}{[}\PY{p}{]}\PY{p}{)}
\end{Verbatim}

    Now, put all about the manifold in one line. The function returns
\textbf{None} if the Dirichlet domain can not be constructed, so the
manifold is \emph{not} hyperbolic:

    \begin{Verbatim}[commandchars=\\\{\}]
{\color{incolor}In [{\color{incolor}5}]:} \PY{k}{def} \PY{n+nf}{getline}\PY{p}{(}\PY{n}{m}\PY{p}{,} \PY{n}{cutoff}\PY{o}{=}\PY{l+m+mf}{3.0}\PY{p}{,} \PY{n}{numcurves}\PY{o}{=}\PY{l+m+mi}{10}\PY{p}{)}\PY{p}{:}
            \PY{k}{try}\PY{p}{:}
                \PY{n}{namelist} \PY{o}{=} \PY{p}{[}\PY{n}{m}\PY{o}{.}\PY{n}{name}\PY{p}{(}\PY{p}{)}\PY{p}{]}
                \PY{n}{thespec} \PY{o}{=} \PY{n}{m}\PY{o}{.}\PY{n}{length\PYZus{}spectrum}\PY{p}{(}\PY{n}{cutoff}\PY{o}{=}\PY{n}{cutoff}\PY{p}{)}
                \PY{n}{lenlist} \PY{o}{=} \PY{n}{mung\PYZus{}complexes}\PY{p}{(}\PY{n}{mung\PYZus{}all\PYZus{}spec}\PY{p}{(}\PY{n}{thespec}\PY{p}{)}\PY{p}{)}
                \PY{n}{vollist} \PY{o}{=} \PY{p}{[}\PY{n+nb}{float}\PY{p}{(}\PY{n}{m}\PY{o}{.}\PY{n}{volume}\PY{p}{(}\PY{p}{)}\PY{p}{)}\PY{p}{]}
                \PY{n}{homo} \PY{o}{=} \PY{n}{m}\PY{o}{.}\PY{n}{homology}\PY{p}{(}\PY{p}{)}
                \PY{n}{qr} \PY{o}{=} \PY{n+nb}{int}\PY{p}{(}\PY{n}{homo}\PY{o}{.}\PY{n}{betti\PYZus{}number}\PY{p}{(}\PY{p}{)}\PY{p}{)}
                \PY{n}{ed} \PY{o}{=} \PY{p}{[}\PY{n+nb}{int}\PY{p}{(}\PY{n}{i}\PY{p}{)} \PY{k}{for} \PY{n}{i} \PY{o+ow}{in} \PY{n}{homo}\PY{o}{.}\PY{n}{elementary\PYZus{}divisors}\PY{p}{(}\PY{p}{)} \PY{k}{if} \PY{n}{i} \PY{o}{\PYZgt{}} \PY{l+m+mi}{0}\PY{p}{]}
                \PY{n}{torsion} \PY{o}{=} \PY{n}{functools}\PY{o}{.}\PY{n}{reduce}\PY{p}{(}\PY{n}{mul}\PY{p}{,} \PY{n}{ed}\PY{p}{,} \PY{l+m+mi}{1}\PY{p}{)}
                \PY{n}{res} \PY{o}{=} \PY{n}{namelist} \PY{o}{+} \PY{n}{lenlist}\PY{p}{[}\PY{p}{:}\PY{n}{numcurves}\PY{p}{]} \PY{o}{+} \PY{n}{vollist} \PY{o}{+} \PY{p}{[}\PY{n}{qr}\PY{p}{,} \PY{n}{torsion}\PY{p}{]}
                
            \PY{k}{except}\PY{p}{:}
                \PY{n}{res} \PY{o}{=} \PY{k+kc}{None}
            \PY{k}{return} \PY{n}{res}
\end{Verbatim}

    For some reason, it's faster to first read in a list, then iterate over
it:

    \begin{Verbatim}[commandchars=\\\{\}]
{\color{incolor}In [{\color{incolor}6}]:} \PY{n}{zoo}\PY{o}{=} \PY{n+nb}{list}\PY{p}{(}\PY{n}{OrientableClosedCensus}\PY{p}{)}
\end{Verbatim}

    Now, read everything in:

    \begin{Verbatim}[commandchars=\\\{\}]
{\color{incolor}In [{\color{incolor}7}]:} \PY{n}{thelines} \PY{o}{=} \PY{p}{[}\PY{n}{getline}\PY{p}{(}\PY{n}{i}\PY{p}{)} \PY{k}{for} \PY{n}{i} \PY{o+ow}{in} \PY{n}{zoo}\PY{p}{]}
\end{Verbatim}

    \begin{Verbatim}[commandchars=\\\{\}]
{\color{incolor}In [{\color{incolor}12}]:} \PY{n}{thelines} \PY{o}{=} \PY{p}{[}\PY{n}{i} \PY{k}{for} \PY{n}{i} \PY{o+ow}{in} \PY{n}{thelines} \PY{k}{if} \PY{n}{i} \PY{o+ow}{is} \PY{o+ow}{not} \PY{k+kc}{None}\PY{p}{]}
\end{Verbatim}

    \begin{Verbatim}[commandchars=\\\{\}]
{\color{incolor}In [{\color{incolor}13}]:} \PY{n}{linedf} \PY{o}{=} \PY{n}{pd}\PY{o}{.}\PY{n}{DataFrame}\PY{o}{.}\PY{n}{from\PYZus{}records}\PY{p}{(}\PY{n}{thelines}\PY{p}{,} \PY{n}{columns} \PY{o}{=} \PY{p}{[}\PY{l+s+s1}{\PYZsq{}}\PY{l+s+s1}{name}\PY{l+s+s1}{\PYZsq{}}\PY{p}{,} \PY{l+s+s1}{\PYZsq{}}\PY{l+s+s1}{a}\PY{l+s+s1}{\PYZsq{}}\PY{p}{,} \PY{l+s+s1}{\PYZsq{}}\PY{l+s+s1}{A}\PY{l+s+s1}{\PYZsq{}}\PY{p}{,} \PY{l+s+s1}{\PYZsq{}}\PY{l+s+s1}{b}\PY{l+s+s1}{\PYZsq{}}\PY{p}{,} \PY{l+s+s1}{\PYZsq{}}\PY{l+s+s1}{B}\PY{l+s+s1}{\PYZsq{}}\PY{p}{,} \PY{l+s+s1}{\PYZsq{}}\PY{l+s+s1}{c}\PY{l+s+s1}{\PYZsq{}}\PY{p}{,} \PY{l+s+s1}{\PYZsq{}}\PY{l+s+s1}{C}\PY{l+s+s1}{\PYZsq{}}\PY{p}{,} \PY{l+s+s1}{\PYZsq{}}\PY{l+s+s1}{d}\PY{l+s+s1}{\PYZsq{}}\PY{p}{,} \PY{l+s+s1}{\PYZsq{}}\PY{l+s+s1}{D}\PY{l+s+s1}{\PYZsq{}}\PY{p}{,} \PY{l+s+s1}{\PYZsq{}}\PY{l+s+s1}{e}\PY{l+s+s1}{\PYZsq{}}\PY{p}{,} \PY{l+s+s1}{\PYZsq{}}\PY{l+s+s1}{E}\PY{l+s+s1}{\PYZsq{}}\PY{p}{,} \PY{l+s+s1}{\PYZsq{}}\PY{l+s+s1}{volume}\PY{l+s+s1}{\PYZsq{}}\PY{p}{,} \PY{l+s+s1}{\PYZsq{}}\PY{l+s+s1}{betti}\PY{l+s+s1}{\PYZsq{}}\PY{p}{,}\PY{l+s+s1}{\PYZsq{}}\PY{l+s+s1}{torsion}\PY{l+s+s1}{\PYZsq{}}\PY{p}{]}\PY{p}{)}
\end{Verbatim}

    \begin{Verbatim}[commandchars=\\\{\}]
{\color{incolor}In [{\color{incolor}14}]:} \PY{n+nb}{len}\PY{p}{(}\PY{n}{zoo}\PY{p}{)}
\end{Verbatim}

\begin{Verbatim}[commandchars=\\\{\}]
{\color{outcolor}Out[{\color{outcolor}14}]:} 11031
\end{Verbatim}
            
    \begin{Verbatim}[commandchars=\\\{\}]
{\color{incolor}In [{\color{incolor}15}]:} \PY{n+nb}{len}\PY{p}{(}\PY{n}{thelines}\PY{p}{)}
\end{Verbatim}

\begin{Verbatim}[commandchars=\\\{\}]
{\color{outcolor}Out[{\color{outcolor}15}]:} 10963
\end{Verbatim}
        \section{First Results}    
    We see that out of the 11031 manifolds, 68 are not hyperbolic.

    The lower case letters columns are of the shortest (5) geodesics, and
the capital letters are of the twisting.

    How about the homology, how is that distributed?

    \begin{Verbatim}[commandchars=\\\{\}]
{\color{incolor}In [{\color{incolor}17}]:} \PY{n}{linedf}\PY{o}{.}\PY{n}{betti}\PY{o}{.}\PY{n}{value\PYZus{}counts}\PY{p}{(}\PY{p}{)}
\end{Verbatim}

\begin{Verbatim}[commandchars=\\\{\}]
{\color{outcolor}Out[{\color{outcolor}17}]:} 0    10836
         1      126
         2        1
         Name: betti, dtype: int64
\end{Verbatim}
            
    We see that 10836 out of the 10963 hyperbolic manifolds (or close to
99\%) are rational homology spheres. This is consistent with the
Dunfield-Thurston model, and also with random fibered manifolds having
\(b_1=1\) with probability approaching \(1.\)

    There are a number of results (notably by Culler and Shalen with
co-authors) on the influence of Betti numbers on volume. Let's see what
we see here:

    \begin{Verbatim}[commandchars=\\\{\}]
{\color{incolor}In [{\color{incolor}18}]:} \PY{n}{linedf}\PY{o}{.}\PY{n}{boxplot}\PY{p}{(}\PY{n}{column}\PY{o}{=}\PY{p}{[}\PY{l+s+s1}{\PYZsq{}}\PY{l+s+s1}{volume}\PY{l+s+s1}{\PYZsq{}}\PY{p}{]}\PY{p}{,} \PY{n}{by}\PY{o}{=}\PY{l+s+s1}{\PYZsq{}}\PY{l+s+s1}{betti}\PY{l+s+s1}{\PYZsq{}}\PY{p}{)}
\end{Verbatim}

\begin{Verbatim}[commandchars=\\\{\}]
{\color{outcolor}Out[{\color{outcolor}18}]:} <matplotlib.axes.\_subplots.AxesSubplot at 0x1d270e45a90>
\end{Verbatim}
            
    \begin{center}
    \adjustimage{max size={0.9\linewidth}{0.9\paperheight}}{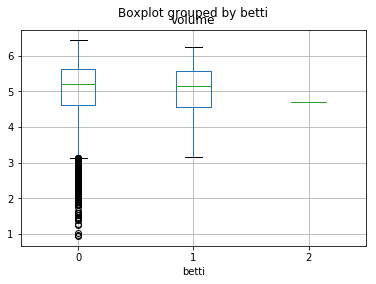}
    \end{center}
    { \hspace*{\fill} \\}
    
    Of course, the number of samples is small in the higher Betti numbers
groups, but it is interesting that the volume \emph{decreases} as the
Betti number increases.

    What about torsion?

    \begin{Verbatim}[commandchars=\\\{\}]
{\color{incolor}In [{\color{incolor}19}]:} \PY{n}{linedf}\PY{o}{.}\PY{n}{torsion}\PY{o}{.}\PY{n}{hist}\PY{p}{(}\PY{n}{bins}\PY{o}{=}\PY{l+m+mi}{30}\PY{p}{)}
\end{Verbatim}

\begin{Verbatim}[commandchars=\\\{\}]
{\color{outcolor}Out[{\color{outcolor}19}]:} <matplotlib.axes.\_subplots.AxesSubplot at 0x1d2713358d0>
\end{Verbatim}
            
    \begin{center}
    \adjustimage{max size={0.9\linewidth}{0.9\paperheight}}{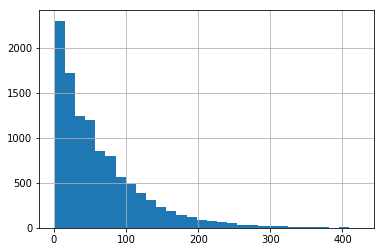}
    \end{center}
    { \hspace*{\fill} \\}
    
    We notice that the highest values of torsion are quite large (given how
small our manifolds are). Let's see how torsion and volume are related.

    \begin{Verbatim}[commandchars=\\\{\}]
{\color{incolor}In [{\color{incolor}20}]:} \PY{n}{sns}\PY{o}{.}\PY{n}{jointplot}\PY{p}{(}\PY{n}{x}\PY{o}{=}\PY{l+s+s1}{\PYZsq{}}\PY{l+s+s1}{torsion}\PY{l+s+s1}{\PYZsq{}}\PY{p}{,} \PY{n}{y}\PY{o}{=}\PY{l+s+s1}{\PYZsq{}}\PY{l+s+s1}{volume}\PY{l+s+s1}{\PYZsq{}}\PY{p}{,} \PY{n}{data}\PY{o}{=}\PY{n}{linedf}\PY{p}{)}
\end{Verbatim}

\begin{Verbatim}[commandchars=\\\{\}]
{\color{outcolor}Out[{\color{outcolor}20}]:} <seaborn.axisgrid.JointGrid at 0x1d271397e10>
\end{Verbatim}
            
    \begin{center}
    \adjustimage{max size={0.9\linewidth}{0.9\paperheight}}{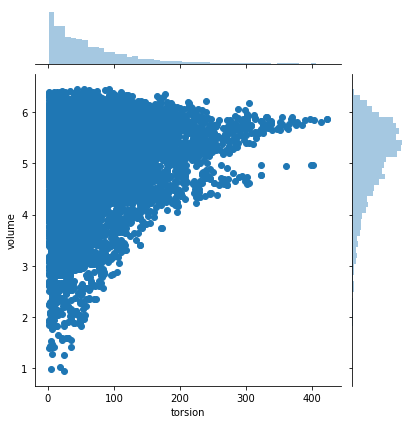}
    \end{center}
    { \hspace*{\fill} \\}
    
    We see that high torsion leads to large volume, though not the other way
around. This is different from the Dunfield-Thurston and random fibered
models, where both volume and log of torsion grow linearly with
complexity. Speaking of log torsion, let's take a look.

    \begin{Verbatim}[commandchars=\\\{\}]
{\color{incolor}In [{\color{incolor}21}]:} \PY{n}{linedf}\PY{p}{[}\PY{l+s+s1}{\PYZsq{}}\PY{l+s+s1}{logtor}\PY{l+s+s1}{\PYZsq{}}\PY{p}{]}\PY{o}{=}\PY{n}{np}\PY{o}{.}\PY{n}{log}\PY{p}{(}\PY{n}{linedf}\PY{o}{.}\PY{n}{torsion}\PY{p}{)}
\end{Verbatim}

    \begin{Verbatim}[commandchars=\\\{\}]
{\color{incolor}In [{\color{incolor}22}]:} \PY{n}{sns}\PY{o}{.}\PY{n}{jointplot}\PY{p}{(}\PY{n}{x}\PY{o}{=}\PY{l+s+s1}{\PYZsq{}}\PY{l+s+s1}{logtor}\PY{l+s+s1}{\PYZsq{}}\PY{p}{,} \PY{n}{y}\PY{o}{=}\PY{l+s+s1}{\PYZsq{}}\PY{l+s+s1}{volume}\PY{l+s+s1}{\PYZsq{}}\PY{p}{,} \PY{n}{data} \PY{o}{=} \PY{n}{linedf}\PY{p}{)}
\end{Verbatim}

\begin{Verbatim}[commandchars=\\\{\}]
{\color{outcolor}Out[{\color{outcolor}22}]:} <seaborn.axisgrid.JointGrid at 0x1d270e45828>
\end{Verbatim}
            
    \begin{center}
    \adjustimage{max size={0.9\linewidth}{0.9\paperheight}}{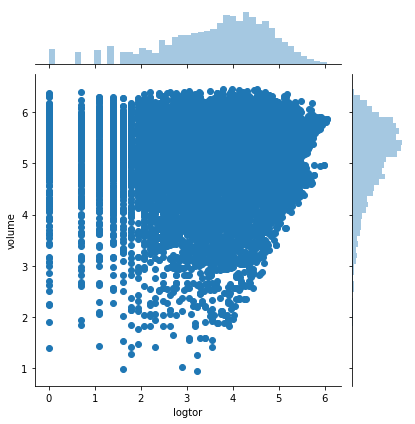}
    \end{center}
    { \hspace*{\fill} \\}
    
    We see that both volume and log torsion are potentially tending to a
gaussian distribution. We also see the linear density cutoff in the
center graph, showing that the random models do have something going for
them. The graph also seems to indicate that
\(V(M) \geq a \log \mbox{tor}(M) + b,\) for some \(a > 0.\)
\section{Pairwise relationships}

    \begin{Verbatim}[commandchars=\\\{\}]
{\color{incolor}In [{\color{incolor}24}]:} \PY{n}{sns}\PY{o}{.}\PY{n}{pairplot}\PY{p}{(}\PY{n}{linedf}\PY{p}{)}
\end{Verbatim}

\begin{Verbatim}[commandchars=\\\{\}]
{\color{outcolor}Out[{\color{outcolor}24}]:} <seaborn.axisgrid.PairGrid at 0x1d2727ee908>
\end{Verbatim}
            
    \begin{center}
    \adjustimage{max size={0.9\linewidth}{0.9\paperheight}}{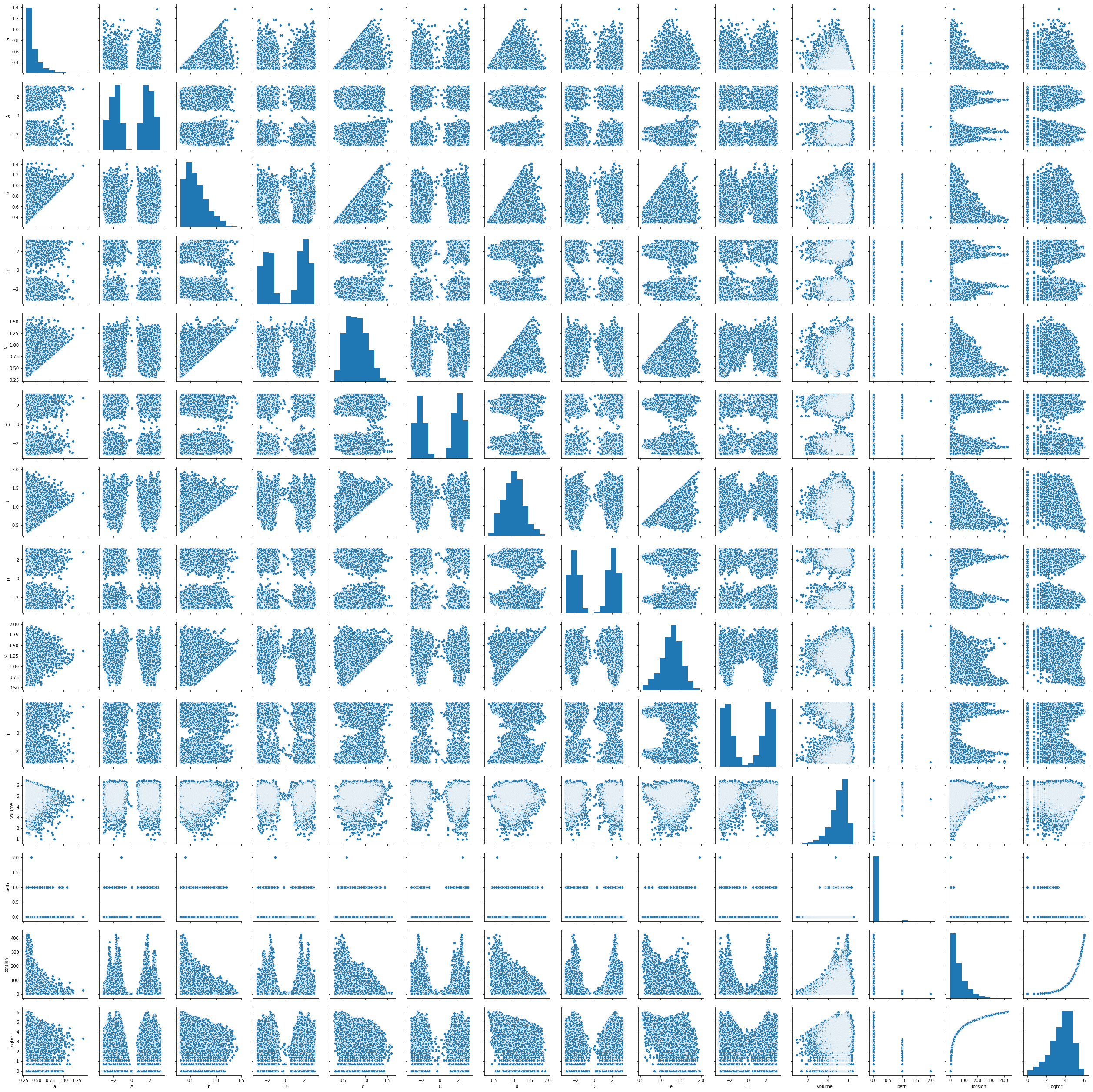}
    \end{center}
    { \hspace*{\fill} \\}
    
    Above, on the diagonal we have the histograms of the various columns and
the off-diagonal cells are the scatter plots of columns against one
another. We see many interesting phenomena.

    \begin{enumerate}
\def\labelenumi{\arabic{enumi}.}
\tightlist
\item
  Notice that the imaginary parts of the complex lengths (the rotations)
  are NOT equidistributed - they are well-separated from 0 (and \(\pi\))
  - they seem to become somewhat less so for longer geodesics.
\item
  The length of the systole (shortest geodesic) is exponentially
  distributed, while
\item
  The lengths of the \(k\)-th shortest geodesics become more and more
  normally distributed as \(k\) becomes large.
\item
  Both the real and the imaginary parts of the lengths seem to be
  uncorrelated (aside from the obvious relation of ordering on the real
  parts).
\item
  Log of torsion seems more-or-less normally distributed.
\end{enumerate}

    The scatter graphs of volume vs the other observables are also
interesting. Let's see if there is any volume between volume and the
systole:

    \begin{Verbatim}[commandchars=\\\{\}]
{\color{incolor}In [{\color{incolor}25}]:} \PY{n}{sns}\PY{o}{.}\PY{n}{regplot}\PY{p}{(}\PY{n}{x}\PY{o}{=}\PY{l+s+s1}{\PYZsq{}}\PY{l+s+s1}{a}\PY{l+s+s1}{\PYZsq{}}\PY{p}{,} \PY{n}{y}\PY{o}{=}\PY{l+s+s1}{\PYZsq{}}\PY{l+s+s1}{volume}\PY{l+s+s1}{\PYZsq{}}\PY{p}{,} \PY{n}{data}\PY{o}{=}\PY{n}{linedf}\PY{p}{)}
\end{Verbatim}

\begin{Verbatim}[commandchars=\\\{\}]
{\color{outcolor}Out[{\color{outcolor}25}]:} <matplotlib.axes.\_subplots.AxesSubplot at 0x1d27af28710>
\end{Verbatim}
            
    \begin{center}
    \adjustimage{max size={0.9\linewidth}{0.9\paperheight}}{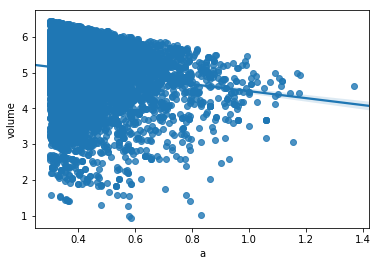}
    \end{center}
    { \hspace*{\fill} \\}
    
    While knowing that the systole is short gives us little information on
the volume, knowing that it is \emph{long} tells us that the volume is
\emph{small}, which is a bit counter-intuitive, at least to this author.
Also, it is pretty that an inequality of the vorm
\(V(M) \leq a s(M) + b,\) where \(s(M)\) is the systole length, and
\(a<0,\) holds.

    What about other geodesics?

    \begin{Verbatim}[commandchars=\\\{\}]
{\color{incolor}In [{\color{incolor}26}]:} \PY{n}{sns}\PY{o}{.}\PY{n}{regplot}\PY{p}{(}\PY{n}{x}\PY{o}{=}\PY{l+s+s1}{\PYZsq{}}\PY{l+s+s1}{c}\PY{l+s+s1}{\PYZsq{}}\PY{p}{,} \PY{n}{y}\PY{o}{=}\PY{l+s+s1}{\PYZsq{}}\PY{l+s+s1}{volume}\PY{l+s+s1}{\PYZsq{}}\PY{p}{,} \PY{n}{data}\PY{o}{=}\PY{n}{linedf}\PY{p}{)}
\end{Verbatim}

\begin{Verbatim}[commandchars=\\\{\}]
{\color{outcolor}Out[{\color{outcolor}26}]:} <matplotlib.axes.\_subplots.AxesSubplot at 0x1d27e1cff98>
\end{Verbatim}
            
    \begin{center}
    \adjustimage{max size={0.9\linewidth}{0.9\paperheight}}{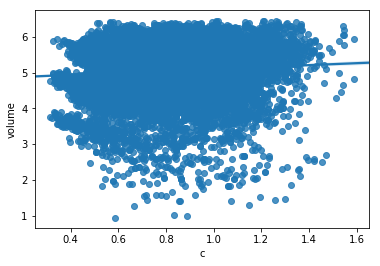}
    \end{center}
    { \hspace*{\fill} \\}
    
    The second longest curve has some slight (but \emph{positive})
predictive power.

    \begin{Verbatim}[commandchars=\\\{\}]
{\color{incolor}In [{\color{incolor}27}]:} \PY{n}{sns}\PY{o}{.}\PY{n}{regplot}\PY{p}{(}\PY{n}{x}\PY{o}{=}\PY{l+s+s1}{\PYZsq{}}\PY{l+s+s1}{d}\PY{l+s+s1}{\PYZsq{}}\PY{p}{,} \PY{n}{y}\PY{o}{=}\PY{l+s+s1}{\PYZsq{}}\PY{l+s+s1}{volume}\PY{l+s+s1}{\PYZsq{}}\PY{p}{,} \PY{n}{data}\PY{o}{=}\PY{n}{linedf}\PY{p}{)}
\end{Verbatim}

\begin{Verbatim}[commandchars=\\\{\}]
{\color{outcolor}Out[{\color{outcolor}27}]:} <matplotlib.axes.\_subplots.AxesSubplot at 0x1d27e1d1b70>
\end{Verbatim}
            
    \begin{center}
    \adjustimage{max size={0.9\linewidth}{0.9\paperheight}}{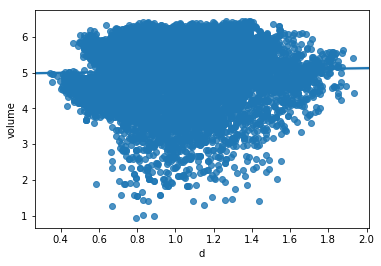}
    \end{center}
    { \hspace*{\fill} \\}
    
    Same with the third\ldots{}

\section{The space of 3-manifolds}
    The question is now, whether we can deduce any of the observables if we
know the others. For example, can we predict the volume from the length
spectrum?

    There is an existence proof, and a construction. For existence, let's
see if there are approximate linear relationships between our many
fields (let's drop the homological invariants for now):

    \begin{Verbatim}[commandchars=\\\{\}]
{\color{incolor}In [{\color{incolor}28}]:} \PY{n}{linedfmin} \PY{o}{=} \PY{n}{linedf}\PY{p}{[}\PY{p}{[}\PY{l+s+s1}{\PYZsq{}}\PY{l+s+s1}{a}\PY{l+s+s1}{\PYZsq{}}\PY{p}{,} \PY{l+s+s1}{\PYZsq{}}\PY{l+s+s1}{A}\PY{l+s+s1}{\PYZsq{}}\PY{p}{,} \PY{l+s+s1}{\PYZsq{}}\PY{l+s+s1}{b}\PY{l+s+s1}{\PYZsq{}}\PY{p}{,} \PY{l+s+s1}{\PYZsq{}}\PY{l+s+s1}{B}\PY{l+s+s1}{\PYZsq{}}\PY{p}{,} \PY{l+s+s1}{\PYZsq{}}\PY{l+s+s1}{c}\PY{l+s+s1}{\PYZsq{}}\PY{p}{,} \PY{l+s+s1}{\PYZsq{}}\PY{l+s+s1}{C}\PY{l+s+s1}{\PYZsq{}}\PY{p}{,} \PY{l+s+s1}{\PYZsq{}}\PY{l+s+s1}{d}\PY{l+s+s1}{\PYZsq{}}\PY{p}{,} \PY{l+s+s1}{\PYZsq{}}\PY{l+s+s1}{D}\PY{l+s+s1}{\PYZsq{}}\PY{p}{,} \PY{l+s+s1}{\PYZsq{}}\PY{l+s+s1}{e}\PY{l+s+s1}{\PYZsq{}}\PY{p}{,} \PY{l+s+s1}{\PYZsq{}}\PY{l+s+s1}{E}\PY{l+s+s1}{\PYZsq{}}\PY{p}{,} \PY{l+s+s1}{\PYZsq{}}\PY{l+s+s1}{volume}\PY{l+s+s1}{\PYZsq{}}\PY{p}{]}\PY{p}{]}
\end{Verbatim}

    \begin{Verbatim}[commandchars=\\\{\}]
{\color{incolor}In [{\color{incolor}29}]:} \PY{n}{u}\PY{p}{,} \PY{n}{s}\PY{p}{,} \PY{n}{v} \PY{o}{=} \PY{n}{np}\PY{o}{.}\PY{n}{linalg}\PY{o}{.}\PY{n}{svd}\PY{p}{(}\PY{n}{linedfmin}\PY{p}{)}
\end{Verbatim}

    \begin{Verbatim}[commandchars=\\\{\}]
{\color{incolor}In [{\color{incolor}30}]:} \PY{n}{s}
\end{Verbatim}

\begin{Verbatim}[commandchars=\\\{\}]
{\color{outcolor}Out[{\color{outcolor}30}]:} array([573.37142695, 257.38960156, 238.87874045, 218.20623087,
                209.10291146, 203.66714934,  45.86448121,  24.94621714,
                 17.98065278,  13.50396573,  10.96093658])
\end{Verbatim}
            
    We see that most of the energy is contained in the top six singular
values, so the space is approximately six dimensional (as a
\emph{linear} space).

    Are the twist parameters redundant somehow?

    \begin{Verbatim}[commandchars=\\\{\}]
{\color{incolor}In [{\color{incolor}51}]:} \PY{n}{linedfmin2} \PY{o}{=} \PY{n}{linedf}\PY{p}{[}\PY{p}{[}\PY{l+s+s1}{\PYZsq{}}\PY{l+s+s1}{a}\PY{l+s+s1}{\PYZsq{}}\PY{p}{,} \PY{l+s+s1}{\PYZsq{}}\PY{l+s+s1}{b}\PY{l+s+s1}{\PYZsq{}}\PY{p}{,}  \PY{l+s+s1}{\PYZsq{}}\PY{l+s+s1}{c}\PY{l+s+s1}{\PYZsq{}}\PY{p}{,}  \PY{l+s+s1}{\PYZsq{}}\PY{l+s+s1}{d}\PY{l+s+s1}{\PYZsq{}}\PY{p}{,}  \PY{l+s+s1}{\PYZsq{}}\PY{l+s+s1}{e}\PY{l+s+s1}{\PYZsq{}}\PY{p}{,} \PY{l+s+s1}{\PYZsq{}}\PY{l+s+s1}{volume}\PY{l+s+s1}{\PYZsq{}}\PY{p}{]}\PY{p}{]}
\end{Verbatim}

    \begin{Verbatim}[commandchars=\\\{\}]
{\color{incolor}In [{\color{incolor}52}]:} \PY{n}{u}\PY{p}{,} \PY{n}{s}\PY{p}{,} \PY{n}{v} \PY{o}{=} \PY{n}{np}\PY{o}{.}\PY{n}{linalg}\PY{o}{.}\PY{n}{svd}\PY{p}{(}\PY{n}{linedfmin2}\PY{p}{)}
\end{Verbatim}

    \begin{Verbatim}[commandchars=\\\{\}]
{\color{incolor}In [{\color{incolor}53}]:} \PY{n}{s}
\end{Verbatim}

\begin{Verbatim}[commandchars=\\\{\}]
{\color{outcolor}Out[{\color{outcolor}53}]:} array([572.95525389,  45.89914386,  24.99847935,  17.9864745 ,
                 13.50661835,  10.96897694])
\end{Verbatim}
            
    Not at all! It looks like the space of volume and the five lengths, only
one or two dimensions are signifcant!

    Let's try to go another way and see if knowing the length spectrum we
can predict the volume. For this we will use boosting - a very effective
machine learning technique.

    \begin{Verbatim}[commandchars=\\\{\}]
{\color{incolor}In [{\color{incolor}69}]:} \PY{n}{featdf} \PY{o}{=} \PY{n}{linedfmin2}\PY{p}{[}\PY{n}{linedfmin2}\PY{o}{.}\PY{n}{columns}\PY{p}{[}\PY{p}{:}\PY{o}{\PYZhy{}}\PY{l+m+mi}{1}\PY{p}{]}\PY{p}{]}
\end{Verbatim}

    \begin{Verbatim}[commandchars=\\\{\}]
{\color{incolor}In [{\color{incolor}70}]:} \PY{n}{targdf} \PY{o}{=} \PY{n}{linedfmin2}\PY{p}{[}\PY{l+s+s1}{\PYZsq{}}\PY{l+s+s1}{volume}\PY{l+s+s1}{\PYZsq{}}\PY{p}{]}
\end{Verbatim}

    \begin{Verbatim}[commandchars=\\\{\}]
{\color{incolor}In [{\color{incolor}71}]:} \PY{n}{X\PYZus{}train}\PY{p}{,} \PY{n}{X\PYZus{}test}\PY{p}{,} \PY{n}{y\PYZus{}train}\PY{p}{,} \PY{n}{y\PYZus{}test} \PY{o}{=} \PY{n}{train\PYZus{}test\PYZus{}split}\PY{p}{(}\PY{n}{featdf}\PY{p}{,} \PY{n}{targdf}\PY{p}{)}
\end{Verbatim}

    \begin{Verbatim}[commandchars=\\\{\}]
{\color{incolor}In [{\color{incolor}72}]:} \PY{n}{clf1} \PY{o}{=} \PY{n}{xg}\PY{o}{.}\PY{n}{XGBRegressor}\PY{p}{(}\PY{p}{)}
\end{Verbatim}

    \begin{Verbatim}[commandchars=\\\{\}]
{\color{incolor}In [{\color{incolor}73}]:} \PY{n}{parameters} \PY{o}{=} \PY{p}{\PYZob{}}\PY{l+s+s1}{\PYZsq{}}\PY{l+s+s1}{objective}\PY{l+s+s1}{\PYZsq{}}\PY{p}{:}\PY{p}{[}\PY{l+s+s1}{\PYZsq{}}\PY{l+s+s1}{reg:linear}\PY{l+s+s1}{\PYZsq{}}\PY{p}{]}\PY{p}{,}\PY{l+s+s1}{\PYZsq{}}\PY{l+s+s1}{learning\PYZus{}rate}\PY{l+s+s1}{\PYZsq{}}\PY{p}{:} \PY{p}{[}\PY{l+m+mf}{0.1}\PY{p}{,} \PY{o}{.}\PY{l+m+mi}{3}\PY{p}{,} \PY{l+m+mf}{0.5}\PY{p}{]}\PY{p}{,} \PY{l+s+s1}{\PYZsq{}}\PY{l+s+s1}{max\PYZus{}depth}\PY{l+s+s1}{\PYZsq{}}\PY{p}{:} \PY{p}{[} \PY{l+m+mi}{5}\PY{p}{,} \PY{l+m+mi}{6}\PY{p}{,} \PY{l+m+mi}{7}\PY{p}{]}\PY{p}{,} \PY{l+s+s1}{\PYZsq{}}\PY{l+s+s1}{min\PYZus{}child\PYZus{}weight}\PY{l+s+s1}{\PYZsq{}}\PY{p}{:} \PY{p}{[}\PY{l+m+mi}{9}\PY{p}{]}\PY{p}{,}\PY{l+s+s1}{\PYZsq{}}\PY{l+s+s1}{silent}\PY{l+s+s1}{\PYZsq{}}\PY{p}{:} \PY{p}{[}\PY{l+m+mi}{1}\PY{p}{]}\PY{p}{,} \PY{l+s+s1}{\PYZsq{}}\PY{l+s+s1}{subsample}\PY{l+s+s1}{\PYZsq{}}\PY{p}{:} \PY{p}{[} \PY{l+m+mf}{0.6}\PY{p}{]}\PY{p}{,} \PY{l+s+s1}{\PYZsq{}}\PY{l+s+s1}{colsample\PYZus{}bytree}\PY{l+s+s1}{\PYZsq{}}\PY{p}{:} \PY{p}{[}  \PY{l+m+mf}{0.8}\PY{p}{]}\PY{p}{,}\PY{l+s+s1}{\PYZsq{}}\PY{l+s+s1}{n\PYZus{}estimators}\PY{l+s+s1}{\PYZsq{}}\PY{p}{:} \PY{p}{[} \PY{l+m+mi}{1000}\PY{p}{]}\PY{p}{\PYZcb{}}
\end{Verbatim}

    \begin{Verbatim}[commandchars=\\\{\}]
{\color{incolor}In [{\color{incolor}74}]:} \PY{n}{clf} \PY{o}{=} \PY{n}{GridSearchCV}\PY{p}{(}\PY{n}{clf1}\PY{p}{,}\PY{n}{parameters}\PY{p}{,} \PY{n}{cv} \PY{o}{=} \PY{l+m+mi}{2}\PY{p}{,} \PY{n}{n\PYZus{}jobs} \PY{o}{=} \PY{l+m+mi}{5}\PY{p}{,} \PY{n}{verbose}\PY{o}{=}\PY{k+kc}{True}\PY{p}{)}
\end{Verbatim}

    \begin{Verbatim}[commandchars=\\\{\}]
{\color{incolor}In [{\color{incolor}75}]:} \PY{n}{clf}\PY{o}{.}\PY{n}{fit}\PY{p}{(}\PY{n}{X\PYZus{}train}\PY{p}{,}\PY{n}{y\PYZus{}train}\PY{p}{)}
\end{Verbatim}

    \begin{Verbatim}[commandchars=\\\{\}]
Fitting 2 folds for each of 9 candidates, totalling 18 fits

    \end{Verbatim}

    \begin{Verbatim}[commandchars=\\\{\}]
[Parallel(n\_jobs=5)]: Using backend LokyBackend with 5 concurrent workers.
[Parallel(n\_jobs=5)]: Done  18 out of  18 | elapsed:   41.9s finished

    \end{Verbatim}

\begin{Verbatim}[commandchars=\\\{\}]
{\color{outcolor}Out[{\color{outcolor}75}]:} GridSearchCV(cv=2, error\_score='raise-deprecating',
                estimator=XGBRegressor(base\_score=0.5, booster='gbtree', colsample\_bylevel=1,
                colsample\_bytree=1, gamma=0, learning\_rate=0.1, max\_delta\_step=0,
                max\_depth=3, min\_child\_weight=1, missing=None, n\_estimators=100,
                n\_jobs=1, nthread=None, objective='reg:linear', random\_state=0,
                reg\_alpha=0, reg\_lambda=1, scale\_pos\_weight=1, seed=None,
                silent=True, subsample=1),
                fit\_params=None, iid='warn', n\_jobs=5,
                param\_grid=\{'objective': ['reg:linear'], 'learning\_rate': [0.1, 0.3, 0.5], 'max\_depth': [5, 6, 7], 'min\_child\_weight': [9], 'silent': [1], 'subsample': [0.6], 'colsample\_bytree': [0.8], 'n\_estimators': [1000]\},
                pre\_dispatch='2*n\_jobs', refit=True, return\_train\_score='warn',
                scoring=None, verbose=True)
\end{Verbatim}
            
    \begin{Verbatim}[commandchars=\\\{\}]
{\color{incolor}In [{\color{incolor}76}]:} \PY{n}{preds} \PY{o}{=} \PY{n}{clf}\PY{o}{.}\PY{n}{predict}\PY{p}{(}\PY{n}{X\PYZus{}test}\PY{p}{)}
\end{Verbatim}

    \begin{Verbatim}[commandchars=\\\{\}]
{\color{incolor}In [{\color{incolor}101}]:} \PY{n}{np}\PY{o}{.}\PY{n}{linalg}\PY{o}{.}\PY{n}{norm}\PY{p}{(}\PY{n}{preds}\PY{o}{\PYZhy{}}\PY{n}{y\PYZus{}test}\PY{p}{)}\PY{o}{/}\PY{n}{np}\PY{o}{.}\PY{n}{sqrt}\PY{p}{(}\PY{n+nb}{len}\PY{p}{(}\PY{n}{y\PYZus{}test}\PY{p}{)}\PY{p}{)}
\end{Verbatim}

\begin{Verbatim}[commandchars=\\\{\}]
{\color{outcolor}Out[{\color{outcolor}101}]:} 0.7417280255021659
\end{Verbatim}
            
    \begin{Verbatim}[commandchars=\\\{\}]
{\color{incolor}In [{\color{incolor} }]:} 
\end{Verbatim}

    \begin{Verbatim}[commandchars=\\\{\}]
{\color{incolor}In [{\color{incolor}43}]:} \PY{n}{y\PYZus{}test}\PY{o}{.}\PY{n}{describe}\PY{p}{(}\PY{p}{)}
\end{Verbatim}

\begin{Verbatim}[commandchars=\\\{\}]
{\color{outcolor}Out[{\color{outcolor}43}]:} count    2741.000000
         mean        5.061871
         std         0.805381
         min         1.529477
         25\%         4.626565
         50\%         5.228348
         75\%         5.657743
         max         6.453448
         Name: volume, dtype: float64
\end{Verbatim}
            
    So we explain about 10\% of the standard deviation, or 20\% of the
variance\ldots{}

    What if we use all the complex length info?

    \begin{Verbatim}[commandchars=\\\{\}]
{\color{incolor}In [{\color{incolor}103}]:} \PY{n}{featdf} \PY{o}{=} \PY{n}{linedfmin}\PY{p}{[}\PY{n}{linedfmin}\PY{o}{.}\PY{n}{columns}\PY{p}{[}\PY{p}{:}\PY{o}{\PYZhy{}}\PY{l+m+mi}{1}\PY{p}{]}\PY{p}{]}
          \PY{n}{targdf} \PY{o}{=} \PY{n}{linedfmin}\PY{p}{[}\PY{l+s+s1}{\PYZsq{}}\PY{l+s+s1}{volume}\PY{l+s+s1}{\PYZsq{}}\PY{p}{]}
          \PY{n}{X\PYZus{}train}\PY{p}{,} \PY{n}{X\PYZus{}test}\PY{p}{,} \PY{n}{y\PYZus{}train}\PY{p}{,} \PY{n}{y\PYZus{}test} \PY{o}{=} \PY{n}{train\PYZus{}test\PYZus{}split}\PY{p}{(}\PY{n}{featdf}\PY{p}{,} \PY{n}{targdf}\PY{p}{)}
          \PY{n}{clf1} \PY{o}{=} \PY{n}{xg}\PY{o}{.}\PY{n}{XGBRegressor}\PY{p}{(}\PY{p}{)}
          \PY{n}{clf} \PY{o}{=} \PY{n}{GridSearchCV}\PY{p}{(}\PY{n}{clf1}\PY{p}{,}\PY{n}{parameters}\PY{p}{,} \PY{n}{cv} \PY{o}{=} \PY{l+m+mi}{2}\PY{p}{,} \PY{n}{n\PYZus{}jobs} \PY{o}{=} \PY{l+m+mi}{5}\PY{p}{,} \PY{n}{verbose}\PY{o}{=}\PY{k+kc}{True}\PY{p}{)}
          \PY{n}{clf}\PY{o}{.}\PY{n}{fit}\PY{p}{(}\PY{n}{X\PYZus{}train}\PY{p}{,}\PY{n}{y\PYZus{}train}\PY{p}{)}
\end{Verbatim}

    \begin{Verbatim}[commandchars=\\\{\}]
Fitting 2 folds for each of 9 candidates, totalling 18 fits

    \end{Verbatim}

    \begin{Verbatim}[commandchars=\\\{\}]
[Parallel(n\_jobs=5)]: Using backend LokyBackend with 5 concurrent workers.
[Parallel(n\_jobs=5)]: Done  18 out of  18 | elapsed:  1.2min finished

    \end{Verbatim}

\begin{Verbatim}[commandchars=\\\{\}]
{\color{outcolor}Out[{\color{outcolor}103}]:} GridSearchCV(cv=2, error\_score='raise-deprecating',
                 estimator=XGBRegressor(base\_score=0.5, booster='gbtree', colsample\_bylevel=1,
                 colsample\_bytree=1, gamma=0, learning\_rate=0.1, max\_delta\_step=0,
                 max\_depth=3, min\_child\_weight=1, missing=None, n\_estimators=100,
                 n\_jobs=1, nthread=None, objective='reg:linear', random\_state=0,
                 reg\_alpha=0, reg\_lambda=1, scale\_pos\_weight=1, seed=None,
                 silent=True, subsample=1),
                 fit\_params=None, iid='warn', n\_jobs=5,
                 param\_grid=\{'objective': ['reg:linear'], 'learning\_rate': [0.1, 0.3, 0.5], 'max\_depth': [5, 6, 7], 'min\_child\_weight': [9], 'silent': [1], 'subsample': [0.6], 'colsample\_bytree': [0.8], 'n\_estimators': [1000]\},
                 pre\_dispatch='2*n\_jobs', refit=True, return\_train\_score='warn',
                 scoring=None, verbose=True)
\end{Verbatim}
            
    \begin{Verbatim}[commandchars=\\\{\}]
{\color{incolor}In [{\color{incolor}104}]:} \PY{n}{preds} \PY{o}{=} \PY{n}{clf}\PY{o}{.}\PY{n}{predict}\PY{p}{(}\PY{n}{X\PYZus{}test}\PY{p}{)}
          \PY{n}{np}\PY{o}{.}\PY{n}{linalg}\PY{o}{.}\PY{n}{norm}\PY{p}{(}\PY{n}{preds}\PY{o}{\PYZhy{}}\PY{n}{y\PYZus{}test}\PY{p}{)}\PY{o}{/}\PY{n}{np}\PY{o}{.}\PY{n}{sqrt}\PY{p}{(}\PY{n+nb}{len}\PY{p}{(}\PY{n}{y\PYZus{}test}\PY{p}{)}\PY{p}{)}
\end{Verbatim}

\begin{Verbatim}[commandchars=\\\{\}]
{\color{outcolor}Out[{\color{outcolor}104}]:} 0.5423190901708408
\end{Verbatim}
            
    So we are doing pretty well (getting about 40\% of the information).

    \begin{Verbatim}[commandchars=\\\{\}]
{\color{incolor}In [{\color{incolor}63}]:} \PY{k+kn}{from} \PY{n+nn}{sklearn}\PY{n+nn}{.}\PY{n+nn}{linear\PYZus{}model} \PY{k}{import} \PY{n}{LinearRegression}
\end{Verbatim}

    \begin{Verbatim}[commandchars=\\\{\}]
{\color{incolor}In [{\color{incolor}78}]:} \PY{n}{reg} \PY{o}{=} \PY{n}{LinearRegression}\PY{p}{(}\PY{p}{)}
\end{Verbatim}

    \begin{Verbatim}[commandchars=\\\{\}]
{\color{incolor}In [{\color{incolor}79}]:} \PY{n}{reg}\PY{o}{.}\PY{n}{fit}\PY{p}{(}\PY{n}{X\PYZus{}train}\PY{p}{,} \PY{n}{y\PYZus{}train}\PY{p}{)}
\end{Verbatim}

\begin{Verbatim}[commandchars=\\\{\}]
{\color{outcolor}Out[{\color{outcolor}79}]:} LinearRegression(copy\_X=True, fit\_intercept=True, n\_jobs=None,
                  normalize=False)
\end{Verbatim}
            
    \begin{Verbatim}[commandchars=\\\{\}]
{\color{incolor}In [{\color{incolor}80}]:} \PY{n}{rpred} \PY{o}{=} \PY{n}{reg}\PY{o}{.}\PY{n}{predict}\PY{p}{(}\PY{n}{X\PYZus{}test}\PY{p}{)}
\end{Verbatim}

    \begin{Verbatim}[commandchars=\\\{\}]
{\color{incolor}In [{\color{incolor}100}]:} \PY{n}{np}\PY{o}{.}\PY{n}{linalg}\PY{o}{.}\PY{n}{norm}\PY{p}{(}\PY{n}{y\PYZus{}test}\PY{o}{\PYZhy{}}\PY{n}{rpred}\PY{p}{)}\PY{o}{/}\PY{n}{np}\PY{o}{.}\PY{n}{sqrt}\PY{p}{(}\PY{n+nb}{len}\PY{p}{(}\PY{n}{y\PYZus{}test}\PY{p}{)}\PY{p}{)}
\end{Verbatim}

\begin{Verbatim}[commandchars=\\\{\}]
{\color{outcolor}Out[{\color{outcolor}100}]:} 0.7917050890601595
\end{Verbatim}
            
    We see that the first three lengths in the length spectrum give us
pretty much all the information!

    \begin{Verbatim}[commandchars=\\\{\}]
{\color{incolor}In [{\color{incolor}84}]:} \PY{n}{newdf} \PY{o}{=} \PY{n}{linedf}\PY{p}{[}\PY{p}{[}\PY{l+s+s1}{\PYZsq{}}\PY{l+s+s1}{logtor}\PY{l+s+s1}{\PYZsq{}}\PY{p}{,} \PY{l+s+s1}{\PYZsq{}}\PY{l+s+s1}{a}\PY{l+s+s1}{\PYZsq{}}\PY{p}{,} \PY{l+s+s1}{\PYZsq{}}\PY{l+s+s1}{volume}\PY{l+s+s1}{\PYZsq{}}\PY{p}{]}\PY{p}{]}
\end{Verbatim}

    \begin{Verbatim}[commandchars=\\\{\}]
{\color{incolor}In [{\color{incolor}85}]:} \PY{n}{u}\PY{p}{,} \PY{n}{s}\PY{p}{,} \PY{n}{v} \PY{o}{=} \PY{n}{np}\PY{o}{.}\PY{n}{linalg}\PY{o}{.}\PY{n}{svd}\PY{p}{(}\PY{n}{newdf}\PY{p}{)}
\end{Verbatim}

    \begin{Verbatim}[commandchars=\\\{\}]
{\color{incolor}In [{\color{incolor}86}]:} \PY{n}{s}
\end{Verbatim}

\begin{Verbatim}[commandchars=\\\{\}]
{\color{outcolor}Out[{\color{outcolor}86}]:} array([659.45727388, 106.79900214,  15.42933463])
\end{Verbatim}
            
    \begin{Verbatim}[commandchars=\\\{\}]
{\color{incolor}In [{\color{incolor}105}]:} \PY{n}{reg2} \PY{o}{=} \PY{n}{LinearRegression}\PY{p}{(}\PY{p}{)}
\end{Verbatim}

    \begin{Verbatim}[commandchars=\\\{\}]
{\color{incolor}In [{\color{incolor}107}]:} \PY{n}{reg2}\PY{o}{.}\PY{n}{fit}\PY{p}{(}\PY{n}{newdf}\PY{p}{[}\PY{p}{[}\PY{l+s+s1}{\PYZsq{}}\PY{l+s+s1}{logtor}\PY{l+s+s1}{\PYZsq{}}\PY{p}{,} \PY{l+s+s1}{\PYZsq{}}\PY{l+s+s1}{volume}\PY{l+s+s1}{\PYZsq{}}\PY{p}{]}\PY{p}{]}\PY{p}{,} \PY{n}{newdf}\PY{o}{.}\PY{n}{a}\PY{p}{)}
\end{Verbatim}

\begin{Verbatim}[commandchars=\\\{\}]
{\color{outcolor}Out[{\color{outcolor}107}]:} LinearRegression(copy\_X=True, fit\_intercept=True, n\_jobs=None,
                   normalize=False)
\end{Verbatim}
            
    \begin{Verbatim}[commandchars=\\\{\}]
{\color{incolor}In [{\color{incolor}108}]:} \PY{n}{reg2}\PY{o}{.}\PY{n}{coef\PYZus{}}
\end{Verbatim}

\begin{Verbatim}[commandchars=\\\{\}]
{\color{outcolor}Out[{\color{outcolor}108}]:} array([-0.02357474, -0.01766166])
\end{Verbatim}
            
    \begin{Verbatim}[commandchars=\\\{\}]
{\color{incolor}In [{\color{incolor} }]:} 
\end{Verbatim}

    % Add a bibliography block to the postdoc
    
    \bibliographystyle{alpha}
    \bibliography{census}
    
    \end{document}